\documentclass{jaca}

\usepackage{mathrsfs}
\usepackage[symbol]{footmisc}
\newcommand{\PP}{{\text{\usefont{U}{dsrom}{m}{n}P}}}
\newcommand{\Pas}{$\PP$-a.s in $\Omega$}
\newcommand{\HS}{\text{HS}}
\newcommand{\eps}{\varepsilon}
\newcommand{\pe}{\psi_{\eps}}
\newcommand{\di}{\displaystyle}
\numberwithin{equation}{section}

\begin{document}

% ====================================================================
% Title, authors, abstract, keywords and AMS codes
% ====================================================================

\title[Well-posedness of a time discretization scheme for a stochastic p-Laplace equation]{Well-posedness of a time discretization scheme for a stochastic p-Laplace equation with Neumann boundary conditions}

\author[C. Bauzet, K. Schmitz, C. Sultan and A. Zimmermann]{Caroline Bauzet, Kerstin Schmitz, C\'edric Sultan\\and Aleksandra Zimmermann}

\begin{abstract}
In this contribution, we are interested in the analysis of a semi-implicit time discretization scheme for the approximation of a parabolic equation driven by multiplicative colored noise involving a $p$-Laplace operator (with $p\geq 2$), nonlinear source terms and subject to Neumann boundary conditions. Using the Minty-Browder theorem, we are able to prove the well-posedness of such a scheme.
\end{abstract}

\KeysAndCodes{Stochastic PDE, It\^o integral, multiplicative Lipschitz noise, Euler-Maruyama scheme, p-Laplace operator, Neumann boundary conditions}{60H15, 35K55, 35K92}

% ====================================================================
% Section 1
% ====================================================================

\section{Introduction}\label{sec1}

\subsection{Statement of the problem and assumptions}

Let $(\Omega,\mathcal{F},(\mathcal{F}_t)_{t\geq 0},\PP)$ be a stochastic basis with the usual assumptions on its filtration $(\mathcal{F}_t)_{t\geq 0}$\footnotemark[3], \footnotetext[3]{\textit{i.e.}, $(\mathcal{F}_t)_{t\geq 0}$ is right continuous and $\mathcal{F}_0$ contains all negligible sets of $\mathcal{F}$.} $T>0$ and $D$ be a bounded open domain of $\Rset^d$ where $d\geq 1$. Having in mind applications in solid mechanics for the future (through the study of a $p$-Laplace perturbation of the stochastic nonlinear parabolic differential inclusion studied in \cite{BBBLV}), we consider in this introductory work the Moreau-Yosida approximation of the maximal monotone operator $\partial \mathcal{I}_{[0,1]}$ (see\textit{, e.g.}, \cite{barbu76,brezis73}), denoted by $(\pe)_{\epsilon>0}$ in the sequel, where for any $\epsilon>0$, $\pe: \Rset \rightarrow \Rset$ is defined for all $v\in\Rset$ by
 \begin{align}\label{SPPR}
	\pe(v)=-\frac{(v)^-}{\eps}+\frac{(v-1)^{+}}{\eps}=
	\begin{cases}
		\hfil \di\frac{v}{\eps}, & \text{if } v\leq 0, \\
		\hfil 0, & \text{if } v\in\cci{0,1}, \\
		\hfil \di\frac{v-1}{\eps}, & \text{otherwise}.
	\end{cases}
\end{align}
Then, for a fixed parameter $\epsilon>0$, we are interested in the following problem:
\begin{equation}\label{eqeps}
\left\{
\begin{alignedat}{2}
du_{\epsilon}+\Big(\Delta_p^N u_{\epsilon}+\pe(u_{\epsilon})\Big)\,dt &= G(u_{\epsilon})\,d\mathcal{W}(t)+\left(\beta(u_{\epsilon})+f\right)\,dt &\qquad & \text{in}\quad \Omega \times \ooi{0,T} \times D \\
u_{\epsilon}(0,\cdot) &= u_0 & & \text{in}\quad \Omega\times D \\
\nabla u_{\epsilon} \cdot \vec{n} &= 0 & & \text{on}\quad \Omega \times \ooi{0,T} \times \partial D
\end{alignedat}
\right.
\end{equation}
where $\Delta_p^N u_{\epsilon}=-\div_x{(|\nabla u_{\epsilon}|^{p-2}\nabla u_{\epsilon})}+|u_{\epsilon}|^{p-2}u_{\epsilon}$ is the $p$-Laplace operator (with $p\geq 2$) associated to homogeneous Neumann boundary conditions, $\mathcal{W}$ is a $(\mathcal{F}_t)_{t\geq0}$-adapted $Q$-Wiener process (for a given operator $Q$ defined in the sequel), $G$, $\beta$, $f$ and $u_0$ are given data, and $\vec{n}$ denotes the unit normal vector to $\partial D$ outward to $D$.

Before stating the considered assumptions on these data, let us precise the stochastic framework. We will call $\mathcal{P}_T$ the predictable $\sigma$-field on $\Omega_T$\footnotemark[2]\footnotetext[2]{$\mathcal{P}_{T}:=\sigma(\{ F_s\times (s,t] \ | \ 0\leq s < t \leq T, \ F_s\in \mathcal{F}_s \} \cup \{F_0\times \{0\} \ | \ F_0\in \mathcal{F}_0 \})$ (see \cite[p. 33]{Liu-Rock}). Then, a mapping defined on $\Omega_T$ with values in a separable Banach space $E$ is predictable if it is $\mathcal{P}_{T}$-measurable.} defined by $\Omega_T = \Omega\times \ooi{0,T}$.
By a slight abuse of notation, for any separable Banach space $X$ and $1\leq r<\infty$,
we will denote the space of predictable, $r$-integrable mappings $u:\Omega_T \rightarrow X$ by $L^r(\Omega_T;X)$\textit{, i.e.}, the measurability of $u\in L^r(\Omega_T;X)$ is understood as the predictable measurability, and in particular $u$ is a random variable with values in $L^r(\ooi{0,T};X)$.
Additionally, we fix a separable Hilbert space $U$ such that $L^2(D)\subset U$, a non-negative symmetric trace class operator $Q:U\rightarrow U$ satisfying $Q^{1/2}(U)=L^2(D)$, and an orthonormal basis $(e_j)_j$ of $U$ made of eigenvectors of $Q$ with corresponding eigenvalues $(\lambda_j)_j\subset \coi{0,\infty}$. At last, the separable Hilbert space of Hilbert-Schmidt operators from $L^2(D)$ to $L^2(D)$ is denoted by $\HS(L^2(D))$. Within this stochastic framework, we consider the following hypotheses on the data:
\begin{itemize}
	\item[$(H_1)$:] $u_0\in L^2\left(\Omega;L^2(D)\right)$ is $\mathcal{F}_0$-measurable and verifies $0\leq u_0(\omega,x) \leq 1$, for almost all $(\omega,x)\in \Omega\times D$.
	\item[$(H_2)$:] $G:L^2(D)\rightarrow \HS\left(L^2(D)\right)$ is such that for any $v\in L^2(D)$ and all $j\in\Nset$, $G(v) \left(Q^{1/2}(e_j)\right) = g_j\left(v\right)$ where $g_j: \Rset \rightarrow \Rset$ is continuous and $\operatorname{supp}g_j \subset \ooi{0,1}$ for all $j\in\Nset$.
	Moreover, it is supposed that there exists $L_g\geq 0$ such that for all $r$, $s\in\Rset$,
	\begin{align*}
		\sum_{j=1}^\infty \left|g_j(r)-g_j(s)\right|^2 \leq L_g |r-s|^2. 
	\end{align*}
	In particular, for all $v,w \in L^2(D)$, we have
	\[\|G(v)-G(w)\|^2_{\operatorname{HS}(L^2(D))}\leq L_g\|v-w\|^2_{L^2(D)}.\]
	\item[$(H_3)$:] $\beta: \Rset\rightarrow \Rset$ is a $L_\beta$-Lipschitz-continuous function (with $L_\beta > 0$) such that, for convenience, $\beta(0)=0$. 
	\item[$(H_4)$:] $f\in L^2(\Omega_T;L^2(D))$.
		\item[$(H_5)$:] There exists $(W_j)_j=\left((W_j(t))_{t\geq 0}\right)_{j\in\Nset^\ast}$ a sequence of independent, real-valued Wiener processes with respect to $(\mathcal{F}_t)_{t\geq 0}$ and such that
		\[\forall t\geq 0, \ \mathcal{W}(t):=\sum_{j=1}^{\infty}\sqrt{\lambda_j} e_j W_j(t)=\sum_{j=1}^{\infty}Q^{1/2}(e_j)W_j(t).\] 
\end{itemize}
In the following, let us write $V=W^{1,p}(D)$ endowed with the norm $||v||_V^p=\int_D \left(|\nabla v|^p+|v|^p\right) dx$ and denote its dual space $(W^{1,p}(D))^{\ast}$ by $V^{\ast}$. We are interested in a solution $u_{\epsilon}$ for Problem (\ref{eqeps}) in the sense of the definition below:
\begin{defn}\label{dups} A stochastic process $u_{\epsilon} \in L^2(\Omega;\mathscr{C}(\ooi{0,T};L^2(D)))\cap L^p(\Omega_T;V)$ is a solution to (\ref{eqeps}) if it satisfies
\begin{align*}
u_{\epsilon}(t)-u_0+\int_0^t \Delta_p^N u_{\epsilon}(s)+\psi_{\eps}(u_{\epsilon}(s))\,ds =\int_0^t G(u_{\epsilon}(s))\,d\mathcal{W}(s)+\int_0^t \beta(u_{\epsilon}(s))+f(s) \,ds
\end{align*}
in $L^2(D)$ for all $t\in\cci{0,T}$, and \Pas, where the stochastic integral is defined as in \cite[Section 4.2]{DPZ14}.
\end{defn}

\subsection{State of the art and aim of the study}

We study the existence, uniqueness and measurability of solutions to a semi-implicit Euler-Maruyama time discretization scheme of a parabolic equation involving a $p$-Laplace operator (with $p\geq 2$), nonlinear source terms and subject to homogeneous Neumann boundary conditions. On the right-hand side of the equation, we consider a multiplicative stochastic forcing given by a stochastic It\^{o} integral with respect to a Hilbert space valued $Q$-Wiener process. The penalization term $(\pe)_{\epsilon>0}$ as well as the homogeneous Neumann boundary condition are motivated by applications in solid mechanics, see\textit{, e.g.}, \cite{FREMOND96} and \cite{LRR23} for more details. More precisely, our study can be seen as a first result towards the numerical analysis of a $p$-Laplace-Allen-Cahn type equation with constraint extending the existing contributions \cite{BBBLV, BSVZ24} to the nonlinear setting. In the deterministic case\textit{, i.e.}, without the presence of the stochastic force term, $p$-Laplace equations with Neumann boundary conditions are well-known. Non exhaustively, let us mention \cite{FK95} and the references therein.

In particular, it is well-known that the Dirichlet $p$-Laplace operator is not coercive on the Sobolev space $W^{1,p}$, and instead, the Neumann $p$-Laplace operator may be considered. In the literature on stochastic PDEs, equations with $p$-Laplace type operators are usually addressed with Dirichlet boundary conditions for different values of $1<p<\infty$. In the case of bounded domains, let us mention \cite{BVWZ13, Gess13, Gess13-2, VWZ16} and the list is far from being complete. A semi-implicit Euler-Maruyama time discretization scheme for pseudomonotone SPDEs has been proposed in \cite{VZ19} in the case of multiplicative noise and in \cite{VZ21} for additive noise.

On unbounded domains, let us mention \cite{LLC14, Krause14, SZi23} concerning $p$-Laplace equations subject to additive or multiplicative noise with homogeneous Dirichlet boundary conditions.

Only a few results on stochastic $p$-Laplace equation with Neumann boundary conditions exist in the literature, see\textit{, e.g.}, \cite{Gess11, CT16}.

To the best of our knowledge, a semi-implicit Euler-Maruyama time discretization sche\-me for Problem (\ref{eqeps}) with Neumann boundary conditions has not yet been proposed in the literature. The main result of this contribution is to prove the well-posedness as stated in Proposition \ref{well-posedness_scheme}. Future work will be devoted to the convergence of the scheme as the time discretization parameter and the regularization parameter $\epsilon$ tend towards $0$ simultaneously, in order to obtain existence of a solution to Problem (\ref{eqeps}) with the maximal monotone operator $\partial \mathcal{I}_{[0,1]}$.

% ====================================================================
% Section 2
% ====================================================================

\section{Semi-implicit time discretization scheme}\label{sec2}
\subsection{Presentation of the scheme and main result}
For $M\in\Nset^\ast$, let $0=t_0<t_1<\ldots<t_M=T$ be an equidistant subdivision of the interval $\cci{0,T}$ and let $\tau:=T/M=t_{n+1}-t_{n}$ for all $n\in\llbracket 0,M-1\rrbracket$. For a fixed $\epsilon>0$, we set $u_{\epsilon,0}=u_0$. Then our semi-implicit Euler-Maruyama scheme is constructed as follows: knowing for a fixed $n\in\llbracket 0,M-1\rrbracket$ a $\mathcal{F}_{t_n}$-measurable random variable with values in $L^2(D)$, denoted by $u_{\epsilon,n}$, we search for a $\mathcal{F}_{t_{n+1}}$-measurable random variable with values in $V$, labelled by $u_{\epsilon,{n+1}}$, such that \Pas \ the following equality holds in $L^2(D)$
\begin{equation}\label{240911_01}
	u_{\epsilon,n+1}-u_{\epsilon,n}+\tau\left(\Delta_p^Nu_{\epsilon,n+1}+\pe(u_{\epsilon,n+1})\right)
	=G(u_{\epsilon,n})\Delta_{n+1} \mathcal{W}+\tau\left(\beta(u_{\epsilon,n+1})+f_n\right),
\end{equation} 
where $\Delta_{n+1} \mathcal{W}:=\mathcal{W}(t_{n+1})-\mathcal{W}(t_{n})$ and $\displaystyle f_n = \frac{1}{\tau} \int_{t_n}^{t_{n+1}} f(s) \, ds$, for all $n\in\llbracket 0,M-1\rrbracket$. 

\begin{rem} Note that as
 $u_{\epsilon,n}$ is assumed to be $\mathcal{F}_{t_n}$-measurable, we have 
	\begin{align*}
	G(u_{\epsilon,n})\Delta_{n+1} \mathcal{W}=	\int_{t_{n}}^{t_{n+1}}G(u_{\epsilon,n}) \ d\mathcal{W} &= \sum_{j=1}^{\infty}G(u_{\epsilon,n})(Q^{1/2}e_j)\left(W_j(t_{n+1})-W_j(t_n)\right)\nonumber \\
		&= \sum_{j=1}^{\infty}g_j(u_{\epsilon,n})\left(W_j(t_{n+1})-W_j(t_n)\right).
	\end{align*}
Since $\Delta_{n+1}\mathcal{W}$ takes values in the Hilbert space $U$ and $G(u_{\epsilon,n})\circ Q^{1/2}$ is a Hilbert-Schmidt operator from $U$ to $L^2(D)$, the last infinite sum is a well-defined element of $L^2\left(\Omega;L^2(D)\right)$.
\end{rem}

The main result of this contribution is:
\begin{prop} \label{well-posedness_scheme} Let us assume that Hypotheses $(H_1)$ to $(H_5)$ are satisfied, consider a fixed parameter $\eps>0$, a fixed $M\in \Nset^\ast$, define $\tau=T/M$, $t_n=n\tau$ ($\forall n\in\llbracket 0,M-1\rrbracket$) and $u_{\epsilon,0}=u_0$. Then, for any given $\mathcal{F}_{t_n}$-measurable random variable $u_{\epsilon,n}$ taking values in $L^2(D)$ (where $n\in\llbracket 0,M-1\rrbracket$), under the assumption $\tau<\frac{1}{L_{\beta}}$, there exists a unique $\mathcal{F}_{t_{n+1}}$-measurable random variable $u_{\epsilon,n+1}$ with values in $V$ solving \eqref{240911_01}.
\end{prop}

\subsection{Proof of Proposition \ref{well-posedness_scheme}}

The proof of Proposition \ref{well-posedness_scheme} is based on Minty-Browder theorem, Pettis measurability theorem and calls on an algebraic inequality from \cite[Lemma 2.1., p.107]{Benedetto1989} stated below in Lemma~\ref{AIB}. By denoting for any $x$, $y$ in $\Rset^d$, the euclidean norm of $x$ by $|x|$, and the associated scalar product of $x$ and $y$ by $x \cdot y$, such a lemma reads as follows:
\begin{lem}[Algebraic inequality] \label{AIB} For any $ p \in \coi{2, +\infty}$ and any $d\in \Nset^\ast$, there exists a constant $C(p)>0$ only depending on $p$ such that
\begin{equation*}
 \forall \zeta, \eta \in \Rset^d,\ \left( |\zeta|^{p-2}\zeta - |\eta|^{p-2}\eta \right) \cdot (\zeta - \eta) \geq C(p)\ |\zeta - \eta|^p.
\end{equation*}
\end{lem}
Let us assume that Hypotheses $(H_1)$ to $(H_5)$ are fulfilled, consider a fixed parameter $\eps>0$, a fixed $M\in \Nset^\ast$ and define $\tau=T/M$ satisfying the condition $\tau<\frac{1}{L_{\beta}}$. Denoting by $(.,.)_{L^2(D)}$ the scalar product in $L^2(D)$, we define the operator 
$A_{\eps,\tau}:V\rightarrow V^{\ast}$ by
	\[\langle A_{\eps,\tau}(u),v\rangle_{V^{\ast},V}:=(u,v)_{L^2(D)} +\tau\int_D \left( |\nabla u|^{p-2}\nabla u\cdot\nabla v +|u|^{p-2}uv+\pe(u)v - \beta(u)v \right) \, dx, \ \ \forall u,v\in V.\]
	Since our scheme (\ref{240911_01}) can be rewritten as
	\[A_{\epsilon, \tau}(u_{\epsilon,n+1}) = u_{\epsilon,n} + G(u_{\epsilon,n})\Delta_{n+1}\mathcal{W} + \tau f_n,\]
our objective is to prove in a first step the invertibility of $A_{\eps,\tau}$ by using the theorem of Minty-Browder in order to get the existence of $u_{\epsilon,n+1}$ in $V$, and to show in a second step its $\mathcal{F}_{t_{n+1}}$-measurability thanks to Pettis measurability theorem.\medskip\\
\textbf{Step 1:} Existence of $u_{\epsilon,n+1}$
\begin{itemize}
\item[$\bullet$] Firstly, we show that $A_{\eps,\tau}$ is coercive.
By the monotonicity of $\pe$, the Lipschitz continuity of $\beta$, the fact that $\pe(0) =\beta(0)= 0$, and the assumption $\tau<\frac{1}{L_{\beta}}$, one gets that for any $u$ in $V$,
\begin{align*}
\langle A_{\eps,\tau}(u), u\rangle_{V^{\ast},V}&= \norm{u}_{L^2(D)}^2 + \tau \int_D \left( |\nabla u|^p + |u|^p + \pe(u)u - \beta(u)u \right) \, dx \\
&\geq (1-\tau L_{\beta}) \norm{u}_{L^2(D)}^2+ \tau \norm{u}_V^p\\
&\geq \tau \norm{u}_V^p,
\end{align*}
and, since $p\geq 2$, one obtains that $\displaystyle \frac{\langle A_{\eps,\tau}(u),u\rangle_{V^{\ast},V}}{\norm{u}_V}
		 \underset{\norm{u}_V \to +\infty}{\to} + \infty$, hence that $A_{\eps,\tau}$ is coercive.
\item[$\bullet$] Secondly, we prove that $A_{\eps,\tau}$ is strictly monotone.
For any $u,v \in V$ such that $u \neq v$, owing again to the monotonicity of $\pe$, the Lipschitz continuity of $\beta$, the assumption $\tau<\frac{1}{L_{\beta}}$, and the application of Lemma \ref{AIB}, we have:
\begin{align*}
&\ \langle A_{\eps,\tau}(u) - A_{\eps,\tau}(v), u-v\rangle_{V^{\ast},V}\\
=& \ \langle A_{\eps,\tau}(u), u-v\rangle_{V^{\ast},V} - \langle A_{\eps,\tau}(v), u-v\rangle_{V^{\ast},V} \\
		= &\ (u-v, u-v)_{L^2(D)} + \tau \int_D \left[\left( \psi_\eps(u) - \psi_\eps(v) \right) (u-v) - \left( \beta(u) - \beta(v) \right) (u-v) \right] \, dx\\
		&\ + \tau \int_D \left[ \left( |\nabla u|^{p-2}\nabla u - |\nabla v|^{p-2}\nabla v \right) \cdot \nabla (u-v) + \left( |u|^{p-2} u - |v|^{p-2} v \right) (u-v)\right]\, dx\\
		\geq& \ (1-\tau L_\beta)\norm{u-v}_{L^2(D)}^2+\tau \Big( C(p) \norm{\nabla(u-v)}_{L^p(D)}^p + C(p) \norm{u-v}_{L^p(D)}^p \Big) \\
\geq& \ \tau\ C(p) \norm{u-v}_V^p,
\end{align*}
which leads to the strict monotonicity of $A_{\eps,\tau}$.
\item[$\bullet$] Thirdly, we show that $A_{\eps,\tau}$ is hemicontinuous. To this end, we prove the stronger property of demicontinuity for $A_{\eps,\tau}$\textit{, i.e.}, that for all $u \in V$, for any $(u_k)_{k\in \Nset} \subset V$ such that $u_k \underset{k \to +\infty}{\to} u$ in $V$, then
\[\forall w \in V, \ \langle A_{\eps,\tau}(u_k), w \rangle_{V^{\ast}, V} \underset{k \to +\infty}{\to} \langle A_{\eps,\tau}(u), w \rangle_{V^{\ast}, V}.\]
We consider $u$ in $V$ and $(u_k)_{k\in \Nset} \subset V$, such that $u_k \underset{k \to +\infty}{\to} u$ in $V$. Then we have in particular $u_k \underset{k \to +\infty}{\to} u$ in $L^p(D)$ and $\nabla u_k \underset{k \to +\infty}{\to} \nabla u$ in $(L^p(D))^d$. Note that, since $\psi_\epsilon$ and $\beta$ are Lipschitz continuous, one gets that for any $w\in V$, as $p\geq 2$,
\begin{equation}\label{dcpebeta}
(u_k,w)_{L^2(D)}+\tau \int_D \big(\pe(u_k)-\beta(u_k)\big)w \, dx\underset{k \to +\infty}{\to}(u,w)_{L^2(D)}+\tau \int_D\big( \pe(u)-\beta(u)\big)w \, dx.
\end{equation}
Let us now focus on the analysis of the demicontinuity of the operator $\Delta_p^N: V\rightarrow V^\ast$ defined for any $u$ in $V$ by $\Delta_p^N u=-\div_x(|\nabla u|^{p-2}\nabla u)+|u|^{p-2}u$. By applying the reverse Lebesgue dominated convergence theorem (see \cite[Th. 4.9, p. 94]{BrezisFA}), there exist $g_0, g_1 \in L^p(D)$ such that, up to a not relabeled subsequence of $(u_k)_{k \in \Nset}$,
\begin{itemize}
\item[(i)] $u_k \underset{k \to +\infty}{\to} u$ a.e. on $D$ and $\nabla u_k \underset{k \to +\infty}{\to} \nabla u$ a.e. on $D$.
\item[(ii)] for all $k \in \Nset$ and a.e. in $D$, $|u_k| \leq g_0$ and $|\nabla u_k| \leq g_1$.
\end{itemize}
Thus, by (i), a.e. in $D$, $|u_k|^{p-2}u_k\underset{k \to +\infty}{\to} |u|^{p-2}u$ and $|\nabla u_k|^{p-2}\nabla u_k\underset{k \to +\infty}{\to} |\nabla u|^{p-2}\nabla u$.
Moreover, as $p \geq 1$, by (ii) one obtains that for all $k \in \Nset$ and a.e. in $D$, \[\big||u_k|^{p-2}u_k\big| \leq g_0^{p-1}\text{ and } \big||\nabla u_k|^{p-2}\nabla u_k\big| \leq g_1^{p-1}.\]
Since $g_0, g_1 \in L^p(D)$, then $g_0^{p-1}, g_1^{p-1} \in L^{p'}(D)$ where $p' = \frac{p}{p-1}$, and therefore, by Lebesgue's dominated convergence theorem, 
\begin{equation}\label{cvforteuk}
|u_k|^{p-2}u_k \underset{k \to +\infty}{\to} |u|^{p-2} u \text{ in }L^{p'}(D)\text{ and }|\nabla u_k|^{p-2}\nabla u_k \underset{k \to +\infty}{\to} |\nabla u|^{p-2}\nabla u\text{ in }\left(L^{p'}(D)\right)^d.
\end{equation}
\textit{A priori}, the convergence of $(|u_k|^{p-2}u_k)_{k \in \Nset}$ and of $(|\nabla u_k|^{p-2}\nabla u_k)_{k \in \Nset}$ may be only for a subsequence. But for any subsequences $(|u_{k_l}|^{p-2}u_{k_l})_{l \in \Nset}$ and $(|\nabla u_{k_l}|^{p-2}\nabla u_{k_l})_{l \in \Nset}$ respectively, we may extract a subsequence $(k_{l_j})_{j \in \Nset}$ such that for any $j \in \Nset$, (i) and (ii) hold true, which is enough to obtain (\ref{cvforteuk}) for the whole sequences $(|u_k|^{p-2}u_k)_{k \in \Nset}$ and $(|\nabla u_k|^{p-2}\nabla u_k)_{k \in \Nset}$.\\
Next, for any $w$ in $V$, making use of the divergence theorem, we get
\begin{align}
\langle \Delta_p^N(u_k), w \rangle_{V^\ast, V}
=\ \ & \int_D - \div_x \left( |\nabla u_k|^{p-2}\nabla u_k \right) w + |u_k|^{p-2}u_k w \, dx\nonumber \\
=\ \ & \int_D |\nabla u_k|^{p-2}\nabla u_k \cdot \nabla w \, dx + \int_D |u_k|^{p-2}u_k w \, dx\nonumber \\
\underset{k \to +\infty}{\to}& \int_D |\nabla u|^{p-2}\nabla u \cdot \nabla w \, dx + \int_D |u|^{p-2} u w \, dx\nonumber \\
=\ \ & \langle \Delta_p^N(u), w \rangle_{V^\ast, V} \label{dcdeltap}.
\end{align}
Finally, (\ref{dcpebeta}) combined with (\ref{dcdeltap}) give us the demicontinuity of $A_{\eps,\tau}$, and afterward its hemicontinuity.
\item[$\bullet$] Fourthly, owing to the coercivity, hemicontinuity, and strict monotonicity of the operator $A_{\eps,\tau}$, one is able to apply the Minty-Browder theorem (see\textit{, e.g.}, \cite[Theorem 2.14, p.38]{Roubicek}), which allows us to affirm that $A_{\eps,\tau}$ is invertible, of inverse operator $A_{\eps,\tau}^{-1}$. Consequently, by setting $u_{\epsilon,0}=u_0$ and $t_n=n\tau$, $\forall n\in\llbracket 0,M-1\rrbracket$, then for any given $\mathcal{F}_{t_n}$-measurable random variable $u_{\epsilon,n}$ taking values in $L^2(D)$ (with $n\in\llbracket 0,M-1\rrbracket$), there exists a unique element $u_{\epsilon,n+1}$ in $V$, such that \Pas,
\[u_{\epsilon,n+1} = A_{\eps,\tau}^{-1}\left(u_{\epsilon,n} + G(u_{\epsilon,n}) \Delta_{n+1}\mathcal{W} + \tau f_n \right).\]
and finally $u_{\epsilon,n+1}$ solves (\ref{240911_01}), \textit{a priori} in $V^\ast$, and \textit{a posteriori} in $L^2(D)$ thanks to the regularity of the data.
\end{itemize}
\textbf{Step 2:} $\mathcal{F}_{t_{n+1}}$-measurability of $u_{\epsilon,n+1}:\Omega\rightarrow V$ \medskip\\
\noindent
In the following, let $A^{-1}_{\tau,\epsilon|L^2(D)}$ be the restriction of $A^{-1}_{\tau,\epsilon}$ to $L^2(D)$. Our objective is to prove the measurability of $A^{-1}_{\tau,\epsilon|L^2(D)}:L^2(D)\rightarrow V$ which will lead us to the $\mathcal{F}_{t_{n+1}}$-measurability of $u_{\epsilon,n+1}:\Omega\rightarrow V$. Indeed, assuming that the above claim holds true, denoting the class of Borel sets of $L^2(D)$ by $\mathcal{B}(L^2(D))$ and the class of Borel sets of $V$ by $\mathcal{B}(V)$, respectively, and setting for any $n\in \llbracket 0,M-1\rrbracket$, $v_n = u_{\epsilon,n} + G(u_{\epsilon,n}) \Delta_{n+1}\mathcal{W} + \tau f_n$, one gets that for all $\mathcal{A}\in \mathcal{B}(V)$,
\begin{equation*}
(u_{\epsilon,n+1})^{-1}(\mathcal{A}) = (A^{-1}_{\tau,\epsilon|L^2(D)} \circ v_n)^{-1}(\mathcal{A}) = v_n^{-1} \left( \big(A_{\tau,\epsilon|L^2(D)}^{-1}\big)^{-1}(\mathcal{A}) \right)\in \mathcal{F}_{t_{n+1}},
\end{equation*}
since $\big(A_{\tau,\epsilon|L^2(D)}^{-1}\big)^{-1}(\mathcal{A}) \in \mathcal{B}(L^2(D))$ and $v_n:\Omega\rightarrow L^2(D)$ is $\mathcal{F}_{t_{n+1}}$-measurable.\smallskip\\
From Pettis measurability theorem (see \cite[Th 1.1.6, p. 5]{HNVW16}), since $V$ is separable, in order to prove the measurability of $A^{-1}_{\tau,\epsilon|L^2(D)}:L^2(D)\rightarrow V$, it is sufficient to show its weak measurability\textit{, i.e.}, that for all $v \in V^{\ast}$, the application $\phi_v:L^2(D)\rightarrow \Rset$ defined by
\[\phi_v(w)=\langle v, A_{\tau,\epsilon|L^2(D)}^{-1}(w) \rangle_{V^{\ast},V}, \ \forall w\in L^2(D),\]
is measurable. To do so, we will prove the weak continuity of $A_{\tau,\epsilon|L^2(D)}^{-1}$ by showing that for every sequence $(h_j)_{j\in\Nset} \subset L^2(D)$ converging strongly to some element $h$ in $L^2(D)$, we have
\[\langle v, A_{\tau,\epsilon|L^2(D)}^{-1}(h_j) \rangle_{V^{\ast},V} \underset{j \to +\infty}{\to} \langle v, A_{\tau,\epsilon|L^2(D)}^{-1}(h) \rangle_{V^{\ast},V}\]
for any fixed $v\in V^{\ast}$. This will give us directly the continuity of $\phi_v$ and subsequently its measurability.\\
Hence, we consider $(h_j)_{j\in\Nset} \subset L^2(D)$ and $h \in L^2(D)$ such that $h_j \underset{j \to +\infty}{\to} h$ in $L^2(D)$. Then for any $j \in \Nset$, setting $u_j=A_{\tau,\epsilon|L^2(D)}^{-1}(h_j)$, by application of $A_{\tau,\epsilon}$ we arrive at the following equality in $V^\ast$:
\begin{equation}\label{191224a}
u_j + \tau \left( \Delta_p^Nu_j + \pe(u_j) - \beta(u_j) \right) = h_j.
\end{equation}
Testing (\ref{191224a}) with $u_j \in V$ yields the following energy estimate for any $\delta>0$
\begin{equation*}
||u_j||_{L^2(D)}^2 + \tau ||u_j||_V^p
%\leq ( h_j, u_j)_{L^2(D)} + \tau ( \beta(u_j), u_j )_{L^2(D)}
\leq \frac{1}{2\delta} ||h_j||_{L^2(D)}^2 + \frac{\delta}{2} ||u_j||_{L^2(D)}^2 + \tau L_\beta ||u_j||_{L^2(D)}^2.
\end{equation*}
By choosing $\delta>0$ such that $1-\frac{\delta}{2}-\tau L_\beta\geq 0$, we get $ ||u_j||_V^p\leq \frac{1}{2\delta\tau} ||h_j||_{L^2(D)}^2$, which proves the boundedness of $(u_j)_{j \in \Nset}$ in $V$. Consequently, there exists $u \in V$ such that, up to a not-relabeled subsequence, $(u_j)_{j \in \Nset}$ converges weakly towards $u$ in $V$, and as $p \geq 2$, by compact embedding of $V$ in $L^p(D)$ we also have, up to a not-relabeled subsequence, strong convergence of $(u_j)_{j \in \Nset}$ towards $u$ in $L^2(D)$. Moreover, by continuity of $\beta$ and $\psi_\epsilon$, we then get the strong convergence in $L^2(D)$ of $(\beta(u_j))_{j\in\Nset}$ and $(\psi_\epsilon(u_j))_{j\in\Nset}$ towards $\beta(u)$ and $\psi_\epsilon(u)$, respectively.\\
Additionally, denoting the conjugate exponent of $p$ by $p'$, since \[\norm{|u_j|^{p-2}u_j}_{L^{p'}(D)} = ||u_j||_{L^p(D)}^{p-1}\text{ and }\norm{|\nabla u_j|^{p-2}\nabla u_j}_{(L^{p'}(D))^d} = \norm{\nabla u_j}_{(L^p(D))^d}^{p-1},\]
one gets that 
$(|u_j|^{p-2}u_j)_{j \in \Nset}$ and $(|\nabla u_j|^{p-2}\nabla u_j)_{j \in \Nset}$ are bounded in $L^{p'}(D)$ and $(L^{p'}(D))^d$, respectively. Consequently, there exist $y \in L^{p'}(D)$ and $\vec{z}\in (L^{p'}(D))^d$ that are weak limits (up to not-relabeled subsequences) respectively of $(|u_j|^{p-2}u_j)_{j\in\Nset}$ in $L^{p'}(D)$ and of $(|\nabla u_j|^{p-2}\nabla u_j)_{j\in\Nset}$ in $(L^{p'}(D))^d$.
Using all these weak and strong convergences in \eqref{191224a}, we obtain that for any $v\in V$,
\begin{equation}\label{191224b}
\int_D uv\, dx + \tau \int_D\big( -\div_x \vec{z} + y + \psi_\epsilon(u) - \beta(u) \big) v\, dx= \int_D hv\, dx, 
\end{equation}
and taking $v=u$ in (\ref{191224b}) provides the following energy equality
\begin{equation}\label{191224c}
\norm{u}_{L^2(D)}^2 = - \tau \langle - \div_x \vec{z} + y, u \rangle_{V^{\ast},V} -\tau (\psi_\epsilon(u), u )_{L^2(D)} + \tau ( \beta(u), u )_{L^2(D)} + ( h, u)_{L^2(D)}.
\end{equation}
In parallel, testing again (\ref{191224a}) by $u_j$ leads to
\begin{equation*}
||u_j||_{L^2(D)}^2 + \tau \langle \Delta^N_p (u_j), u_j \rangle_{V^{\ast},V} + \tau ( \psi_\epsilon(u_j), u_j)_{L^2(D)} -\tau (\beta(u_j), u_j)_{L^2(D)} = ( h_j, u_j )_{L^2(D)},
\end{equation*}
and passing to the superior limit with respect to $j$, one arrives at
\begin{equation}\label{191224d}
\norm{u}_{L^2(D)}^2 = - \tau \limsup_{j \to +\infty} \langle \Delta^N_p (u_j), u_j \rangle_{V^{\ast},V} - \tau ( \psi_\epsilon(u), u )_{L^2(D)} + \tau (\beta(u), u )_{L^2(D)} + ( h, u )_{L^2(D)}.
\end{equation}
From (\ref{191224c}) and (\ref{191224d}), one deduces
\begin{equation*}
\limsup_{j \to +\infty} \langle \Delta^N_p (u_j), u_j \rangle_{V^{\ast},V} = \langle - \div_x \vec{z} + y, u \rangle_{V^{\ast},V}
\end{equation*}
Therefore, by Minty's monotonicity argument, we can affirm that $\Delta^N_p (u) = - \div_x \vec{z} + y$.
Back into (\ref{191224b}), we then obtain the following equality in $V^\ast$
\begin{equation*}
u + \tau \left( \Delta^N_p (u) + \psi_\epsilon(u) - \beta(u) \right) = h, 
\end{equation*}
which leads to $u = A_{\tau,\epsilon|L^2(D)}^{-1}(h)$.\\
At last, fixing $v \in V^{\ast}$, owing to the strong convergence of $(u_j)_{j\in\Nset}=(A_{\tau,\epsilon|L^2(D)}^{-1}(h_j))_{j\in\Nset}$ towards $u$ in $V$, we may now conclude that
\begin{equation*}
\lim_{j \to +\infty} \langle v, A_{\tau,\epsilon|L^2(D)}^{-1}(h_j) \rangle_{V^{\ast},V}
= \lim_{j \to +\infty} \langle v, u_j \rangle_{V^{\ast},V}
= \langle v, u \rangle_{V^{\ast},V}
= \langle v, A_{\tau,\epsilon|L^2(D)}^{-1}(h) \rangle_{V^{\ast},V}
\end{equation*}
and we get back the weak continuity of the operator $A_{\tau,\epsilon|L^2(D)}^{-1}:L^2(D)\rightarrow V$, which completes the proof of Proposition \ref{well-posedness_scheme}.

% ====================================================================
% Acknowledgements
% ====================================================================

\ack The authors would like to thank the German Research Foundation, the Programmes for Project-Related Personal Exchange Procope and Procope Plus and the Procope Mobility Program DEU-22-0004 LG1 for financial support.

% ====================================================================
% References
% ====================================================================

% 1. Make your bib file and write its name below.
% 2. Run pdfLaTeX (twice) on the tex file.
% 3. Run BibTeX on the aux file.
% 4. Run LaTeX again a couple of times to get the correct cross-references.
%
% In many LaTeX editors, the above loop can be done just by pressing a couple of buttons.

\bibliographystyle{acmurl}
\bibliography{bibBSSZ}

% ====================================================================
% Address
% ====================================================================

\begin{address}
	Caroline Bauzet and C\'edric Sultan \\
	Aix Marseille Univ, CNRS, Centrale Med, LMA, Marseille, France \\
	4 Impasse Nikola Tesla \\
	13453 Marseille CEDEX 13, France \\
	\texttt{caroline.bauzet@univ-amu.fr} and \texttt{cedric.sultan@univ-amu.fr}
\end{address}

\begin{address}
	Kerstin Schmitz and Aleksandra Zimmermann \\
	TU Clausthal, Institut für Mathematik, Clausthal-Zellerfeld, Germany \\
	Erzstraße 1 \\
	38678 Clausthal-Zellerfeld, Germany \\
	\texttt{kerstin.schmitz@tu-clausthal.de} and \texttt{aleksandra.zimmermann@tu-clausthal.de}
\end{address}

\end{document}